\documentclass[11pt]{article}
	
	\newcommand{\blind}{0}
	
    \makeatletter
    \renewcommand\section{\@startsection {section}{1}{\z@}%
                                       {-3.5ex \@plus -1ex \@minus -.2ex}%
                                       {2.3ex \@plus.2ex}%
                                       {\normalfont\fontfamily{phv}\fontsize{16}{19}\bfseries}}
    \renewcommand\subsection{\@startsection{subsection}{2}{\z@}%
                                         {-3.25ex\@plus -1ex \@minus -.2ex}%
                                         {1.5ex \@plus .2ex}%
                                         {\normalfont\fontfamily{phv}\fontsize{14}{17}\bfseries}}
    \renewcommand\subsubsection{\@startsection{subsubsection}{3}{\z@}%
                                        {-3.25ex\@plus -1ex \@minus -.2ex}%
                                         {1.5ex \@plus .2ex}%
                                         {\normalfont\normalsize\fontfamily{phv}\fontsize{14}{17}\selectfont}}
    \makeatother
	
	\usepackage{amsmath}
	\usepackage{graphicx}
	\usepackage{subcaption}
	\usepackage{multirow}
	\usepackage{enumerate}
	\usepackage[margin=1in]{geometry}
	\usepackage{natbib} 
	\usepackage{url} 
	\usepackage[font=small,skip=0pt]{caption}
	\usepackage{wrapfig}
	\usepackage{array,float}

	\usepackage{graphicx}
    \graphicspath{ {./images/} }
    \usepackage{float}
    \usepackage{caption}
    \usepackage{booktabs}

		\usepackage{amsfonts, amsthm, latexsym, amssymb}
		\usepackage{lineno}
		\usepackage{xcolor}
		\usepackage{soul}
	
\begin{document}
		
		\def\spacingset#1{\renewcommand{\baselinestretch}%
			{#1}\small\normalsize} \spacingset{1}
		
		\if0\blind
		{
			\title{\bf Improving Access to Housing and Supportive Services for Runaway and Homeless Youth: Reducing Vulnerability to Human Trafficking in New York City}
			\author{Yaren Bilge Kaya$^1$, Kayse Lee Maass$^1$, Geri L. Dimas$^2$, Renata Konrad$^3$,\\
			 Andrew C. Trapp$^{2,3}$, Meredith Dank$^4$ \\
			 \\
$^1$ Mechanical and Industrial Engineering, Northeastern University, Boston, MA \\
$^2$ Data Science Program, Worcester Polytechnic Institute, Worcester, MA \\
$^3$ WPI Business School, Worcester Polytechnic Institute, Worcester, MA \\
$^4$ New York University, New York City, NY}
			\date{}
			\maketitle
		} \fi
		
		\if1\blind
		{

            \title{\bf \emph{Improving Access to Housing and Supportive Services for Runaway and Homeless Youth: Reducing Vulnerability to Human Trafficking in New York City}}
			\author{Author information is purposely removed for double-blind review}
			
\bigskip
			\bigskip
			\bigskip
			\begin{center}
				{\LARGE\bf \emph{Improving Access to Housing and Supportive Services for Runaway and Homeless Youth: Reducing Vulnerability to Human Trafficking in New York City}}
			\end{center}
			\medskip
		} \fi
		\bigskip
		
		\vspace{-12mm}
	\begin{abstract}
\noindent Recent estimates indicate that there are over 1 million runaway and homeless youth and young adults (RHY) in the United States (US). Exposure to trauma, violence, and substance abuse, coupled with a lack of community support services, puts homeless youth at high risk of being exploited and trafficked. Although access to safe housing and supportive services such as physical and mental healthcare is an effective response to youth’s vulnerability towards being trafficked, the number of youth experiencing homelessness exceeds the capacity of available housing resources in most US communities. We undertake a RHY-informed, systematic, and data-driven approach to project the collective capacity required by service providers to adequately meet the needs of RHY in New York City, including those most at risk of being trafficked. Our approach involves an integer linear programming model that extends the multiple multidimensional knapsack problem and is informed by partnerships with key stakeholders. The mathematical model allows for time-dependent allocation and capacity expansion, while incorporating stochastic youth arrivals and length of stays, services provided in a periodic fashion, and service delivery time windows. Our RHY and service provider-centered approach is an important step toward meeting the actual, rather than presumed, survival needs of vulnerable youth.
	\end{abstract}
			
	\noindent%

	\spacingset{1.5} 
	
\section{Introduction} \label{s:intro}
\vspace{-3mm} 
\noindent Human trafficking (HT) is the criminal commercial exchange and exploitation of humans for monetary gain or benefit – a globally prevalent, violation of human rights \citep{Gajic-2007}. HT is rarely a one-time event, but rather a process that can be conceptualized as a series of event-related stages over time during which various risks and intervention opportunities may arise: recruitment; possible transit; exploitation; and re-integration \citep{Zimmerman-2011}. A key and largely understudied means to disrupt trafficking networks is through reducing an individual’s vulnerability during any of the four stages. \\
\indent Runaway and homeless youth and young adults (RHY) are particularly vulnerable to exploitation and trafficking \citep{Wright-2021, Hogan-2020}
. Several factors are known to raise the risk of trafficking for this population, such as a history of physical abuse, emotional neglect, and low self-esteem \citep{Hannan-2017}. Such circumstances, coupled with a lack of community support, put RHY at a high risk of HT. One study on child sexual exploitation found that most victims in the study sample had experienced homelessness or persistent housing instability \citep{Curtis-2008}. Labor traffickers are also known to prey on vulnerabilities by providing housing along with the employment opportunity, making it more difficult for youth to leave the trafficking situation because they have nowhere else to go \citep{Bigelsen-2013}. When youth leave a trafficking situation, they need to have near-immediate access to stable shelter and other support services, as well as financial resources to ensure they are not re-victimized and are able to  recover \citep{Duncan-2019}
. In a recent study with individuals who were able to exit child sexual exploitation, the survivors interviewed confirmed the need for social services to provide ongoing safety and basic needs, such as shelter, to fully exit trafficking and exploitative experiences \citep{Bruhns-2018}. For all these reasons, access to appropriate housing can significantly decrease the likelihood of youth experiencing sexual or labor exploitation and trafficking \citep{Davy-2015, Potocky-2010}.  \\
\indent Trafficking prevention and rehabilitative services have been shown to be effective in disrupting trafficking activity by decreasing vulnerability and associated risk factors. Of the varied prevention and rehabilitative services from which those at-risk of trafficking and HT survivors may benefit, access to safe housing is widely agreed to be the most pressing need \citep{Dank-2015, Clawson-2006}
. Yet, while the exact numbers of available beds for RHY in the US is unknown, it is widely acknowledged that demand greatly exceeds supply \citep{Clawson-2009}.\\ 
\indent As a trafficking intervention, shelter provision extends beyond the mere supplying of beds. Shelters are linked to a dynamic and shifting landscape of networked support services such as medical treatment, psycho-social care, education, life-skills training, and legal advocacy that aids in a holistic approach to rehabilitation and trafficking network disruption \citep{Ide-2019}. The ability to engage at-risk populations in shelter services is a critical component of service provision. Programs that provide at-risk youth and trafficking survivors with options and the ability to make their own choices regarding their care further reduce vulnerability \citep{Hopper-2010}. Yet, shelters are unique in the varied services they offer, as well as inclusion and exclusion criteria that dictate who may receive services at the shelter \citep{Clawson-2007}. These restrictions may depend on a shelter's funding sources, state law, or federal policy.  It is thus important to place youth in shelters that can meet the unique profile of needs of at-risk youth and HT survivors.\\
\indent This study is motivated by the appreciation that shelter and associated services disrupt trafficking activity by decreasing vulnerability for those at-risk of trafficking, including HT survivors. We seek to alleviate current capacity limitations and improve access to housing and support services for homeless and unstably housed youth and young adults as a mechanism to reduce the supply of potential HT victims. The study population of this study are homeless, unstably housed and at-risk youth in New York City (NYC) ages 16-24, we note that NYC has the largest such population in the US \citep{Morton-2019}.\\
\indent Informed by partnerships and interviews with RHY who have experienced homelessness, shelter service providers, and the NYC Mayor’s Office, we develop a mixed-integer linear program (MILP) to project the cost-minimizing shelter and service capacity expansion that meets the collective needs of youth. Our RHY-centered approach allows us to provide a much clearer picture of the actual, rather than presumed, needs of homeless youth. To the best of our knowledge, our study is the first to (i) incorporate primary data collection to estimate the support services and resources available to RHY, and (ii) develop an optimization model to determine the cost-minimizing approach to expand capacity under stochastic demand patterns including RHY arrivals, stay durations, and varying service frequencies and intensities. The presented model could be applied to improve access to other public services such as (i) non-profit organizations that provide food assistance, (ii) legal aid services, and (iii) public healthcare services: especially for patients with chronic diseases. \\
\indent The remainder of the paper is structured as follows. We provide a review of the literature in Section~\ref{s:LitSearch}; explain our data collection process in Section~\ref{s:Data}; propose a mixed-integer linear programming framework in Section~\ref{s:model}; present our computational setup, and discuss results, analysis and insights from our computational experiments in Section~\ref{s:experiments}. We conclude by summarizing our contributions, limitations and future research directions in Section~\ref{s:conclusion}.
\vspace{-4mm}
\section{Literature Review} \label{s:LitSearch}
Operations Research (OR) techniques are applicable to challenges in HT, and engineering systems analysis can produce important insights concerning HT operations \citep{Konrad-2017,Caulkins-2019}. Analytical techniques employed in disaster preparedness may be adapted to allocate limited funds and resources for anti-trafficking operations. However, such techniques likely need to be extended to be  effective as trafficking is an ongoing, rather than a single disaster event \citep{He-2016}. Contextual similarities exist between allocating funds for anti-trafficking services and models addressing fund allocation for hunger relief \citep{Orgut-2016}. Similarly, capacity allocation and location models such as \cite{Kilci-2015} and \cite{Li-2012} could be adapted to increase the efficiency and effectiveness of public services provided to HT victims and populations  vulnerable to HT. However, simply adapting existing OR approaches to address HT is insufficient and irresponsible -- the intricacies of individual agency must be incorporated for successful implementation of OR techniques \citep{Konrad-2022}. Our work draws upon, and contributes to, two primary areas: (i) OR for disrupting the supply side of HT networks and (ii) effective and appropriate allocation of scarce resources via capacity allocation and expansion.  
 \vspace*{-4mm}
 \subsection{\emph{OR, Human Trafficking and Homelessness}}
Research related to OR and analytics efforts to reduce homelessness and the supply of sex and labor trafficking victims has only recently received attention and largely focus on three broad areas: the scope of the problem \citep{Kosmas-2020, Brelsford-2018}, frameworks for addressing the crisis \citep{Tezcan-2020,Taylor-2018}
, and appropriate allocation of scarce resources to combat HT and homelessness \citep{Chan-2018, Petry-2021}. A second body of OR literature focuses on detection of hidden HT victims \citep{Keskin-2021, Kapoor-2017}
, movement patterns of covert traffickers \citep{Mcdonald-2021, Yao-2021}
and long-term intervention approaches to prevent trafficking or re-trafficking by increasing the efficiency and effectiveness of public services \citep{Amadasun-2020}, as well as improving access to supportive and rehabilitative services \citep{Azizi-2018}. However, little attention has been paid to using OR techniques to examine the problem of allocating scarce resources for anti-trafficking efforts, with the exceptions of \cite{Konrad-2019} and \cite{Maass-2020}.

\indent To increase trafficking awareness among at-risk populations, \cite{Konrad-2019} proposes a resource allocation model for trafficking prevention programs in Nepal, aiding decision makers in evaluating how to allocate limited funds in the context of trafficking awareness. To improve stabilization and the eventual reintegration of trafficked persons, \cite{Maass-2020} present an optimization model that allocates a budget for locating residential shelters in a manner that maximizes a measure of societal impact. The model evaluates the trade-off in the cost of opening and operating shelters in each location with the health benefits, labor productivity gained, and reduction in criminal justice costs. Although both \cite{Maass-2020} and \cite{ Konrad-2019} investigate HT-related interventions, they emphasize the broader system and the optimal budget allocation rather than the effective capacity allocation of supportive and rehabilitative services provided to vulnerable populations.

\indent Although not explicitly focused on reducing vulnerability to HT, other OR literature addresses the importance of reducing homelessness. \cite{Chan-2018} propose an Artificial Intelligence-based decision maker to support the current long-term housing assignment process of RHY. They consider the resource-constrained assignment of various types of housing programs available to youth with different needs. In another study, \cite{Azizi-2018} investigate how to prioritize heterogeneous homeless youth on a waiting list for scarce housing resources of different types. They use mixed-integer programming and machine learning to present a policy that trades off efficiency, fairness, and interpretability. Both \cite{Chan-2018} and \cite{Azizi-2018} focus on making the housing assignment process of homeless youth more efficient in the presence of scarce resources. Both sets of authors propose  approaches that prioritize certain groups of youth, rather than expanding housing resources to meet a greater number of homeless youth's needs. Moreover, neither of these studies consider that youth may need particular resources at varied levels at different time periods, nor do they incorporate a time component in their models to reflect the dynamic demands seen in practice. Additionally, both studies assume that there is a centralized list of youth who are waiting to be placed into housing and are available to be assigned simultaneously. It has been well documented that for emergency and transitional independent living shelters, youth often tend to self-select the shelter they visit, possibly due to location, convenience, or practical necessity. The aforementioned assumption may therefore limit their study population to youth in search of long-term housing. 
\vspace*{-4mm}
\subsection{\emph{Capacity Expansion and Allocation Models}}
The problem of planning capacity expansion in facilities that provide different products and services has received considerable attention, and many mathematical programming formulations have been proposed. Earlier studies on this topic mainly focused on expanding electric utility capacity, modeling the problem as a linear program \citep{Williamson-1966, Sherali-1982} 
. A number of analytics-based methods have been applied to capacity expansion in healthcare settings. \cite{Lovejoy-2002} investigate how a hospital can best invest in operating room capacity to provide high-quality service while protecting its profitability by considering the trade-off among three performance criteria: wait time, scheduled procedure start-time reliability, and hospital profits. In another study, \cite{Akcali-2006} present a network flow model that incorporates facility performance and budget constraints to determine optimal hospital bed capacity over a finite planning horizon. Similarly, \cite{Woodall-2013} combine simulation and optimization to improve patient flow within the Duke Cancer Institute. These studies are representative of many healthcare planning and scheduling models while featuring an objective function  centered on efficiency or patient satisfaction, such as cost minimization, revenue maximization, and waiting time minimization.\\
\indent While such studies were instrumental in pioneering capacity expansion and allocation models within healthcare systems, they lack a number of critical elements such as the impact on individual and aggregate health outcomes. To fill this gap, \cite{Deo-2013} develop an integrated capacity allocation model that incorporates clinical (disease progression) and operational (capacity constraint) aspects for chronic disease treatment, and investigates how operational decisions can improve population-level health outcomes. As our study examines the chronic problem of youth homelessness in the context of trafficking, we use a similar approach to form our objective function, integrating both system efficiency and youth preference.\\
\indent Against the backdrop of capacity allocation literature, our decision to use capacity expansion models was motivated by an opportunity to consider service delivery time windows and client preferences in a novel manner. Some of the first research efforts to consider service delivery time windows for clients include \cite{Kolen-1987} and \cite{Desrochers-1992}. Both studies focused on routing limited capacity vehicles while minimizing their total distances to fulfill known demands. While these studies were an important step in the time window literature, they are unable to capture the stochastic nature of client demands and preferences. Demand stochasticity is a necessary component of the real-world modeling conditions and is included in more recent studies \citep{Gocgun-2013, Patrick-2008, Jalilvand-2021}. The work of \cite{Gocgun-2013} employs a Markov Decision Process perspective to assign randomly arriving chemotherapy patients to future appointment dates within clinically established time windows. While the methods of \cite{Gocgun-2013} were effective to capture the randomness of demands, none of the existing studies consider varied service needs within different time windows. Furthermore, existing approaches do not consider preferences of the client. In contrast, we (i) introduce the time window concept into our capacity expansion model; (ii)  consider the unique preferences of youth, and (iii) incorporate a variety of service needs within different time windows.\\
\indent Building homeless shelter and service workforce capacity is important to strengthen the ability to deliver effective services to homeless and vulnerable populations \citep{Mullen-2010}. To the best of our knowledge, no OR-based study exists that focuses on building  shelter and service capacity for the homeless as a means to decrease vulnerability to trafficking. To address this gap we project the cost-minimizing capacity to deploy to organizations providing housing and supportive services to RHY youth while considering: (i) multiple organizations that provide multiple capacitated services, (ii) organizations that only serve certain demographics, (iii) stochastic youth arrivals and stay durations, (iv) varying service frequencies and intensities, (v) service delivery time windows, (vi) periodic and non-periodic service provision, and (vii) youth abandonment.
\vspace*{-4mm}
\section{Data and Community Partners} \label{s:Data}
This section presents a high-level overview of the data acquisition process. New York City (NYC) has the highest rate of homelessness in the United States \citep{USDHUD}, with an estimated 14,946 homeless children and 18,370 single adults sleeping in the NYC municipal shelter system on any given night \citep{CHY}. Our study population is at-risk runaway, homeless, and unstably housed youth and young adults ages 16-24 in NYC and we consider service providers as the non-profit organizations that provide shelter services to this population. For brevity, we will refer to these populations as RHY and RHY organizations. To better understand existing resources and capacities of RHY organizations, service needs and preferences of the youth, as well as the practicality  and feasibility of implementation of any model results, we collected data and feedback from multiple stakeholders: RHY organizations, RHY, the New York City Mayor’s Office, and the New York Coalition for Homeless Youth. Regular meetings were held with the New York City Mayor's Office and the New York Coalition for Homeless Youth to obtain feedback and suggestions for improvements throughout the data collection and model building phases. These data were used to create provision and need profiles of service providers and RHY in NYC, respectively. 
\vspace*{-4mm}
\subsection{\emph{RHY Organization Profiles}} \label{s:Provider}
We conducted structured interviews and surveys with five RHY organizations that fund, support, and provide RHY services to assess existing capacities, resources, and nature of different services in NYC. These interviews and surveys revealed the: (i) demographics of the youth served; (ii) types of services offered by each organization; (iii) amount of resources available for different services; (vi) average length of stay of youth, and (v) services outsourced through referrals to other RHY organizations. The RHY organizations we interviewed provide different types of programs to RHY such as crisis/emergency, transitional independent living (TIL), and long-term housing. We focus on TIL programs as they offer housing and support services to youth while youth work towards establishing independence. In particular, TIL programs aim to help youth gain practice in skills such as education, housing, employment, recreation, health, and safety, all of which promote self-sufficiency and independence \citep{Naccarato-2008}. Although we focus on TIL programs, the model we present can be expanded to incorporate other housing services as well.

\indent The average length of stay (LOS) of RHY in TIL is highly dependent on the RHY organization, program, and youth themselves; in our modeling, the LOS of youth follows a normal distribution. Additionally, reasons such as safety issues, mental health problems, strict RHY organization rules, finding stable housing or reuniting with family might cause RHY to abandon the RHY organization earlier than expected. Thus, following the stakeholder recommendations, we assume that a portion of RHY have shorter average length of stay due to abandonment.

\indent COVID-19 protocols and restrictions put immense operating pressure on RHY organizations, limiting our data collection.  Thus, we focused on eight RHY organizations funded by a major funding agency that we interviewed, and supplemented our primary data (interview and survey responses) with publicly available data. Through qualitative examination of the data obtained across these sources, we assigned each organization a profile of characteristics. In total, we created five TIL organization profiles (services provided and the demographics served) using the primary data sources, and three organization profiles using the secondary data sources. These RHY organization demographics and service profiles used in our model can be found in Table 1 and 3 of the Appendix. We also include an \textit{incompatibility set} to represent a hypothetical RHY organization capable of serving any individual regardless of their demographic. This hypothetical organization is a placeholder for any RHY who cannot be assigned to an existing organization due to demographic mismatch. In practice, it can be used to assess shortcomings in the existing resources to serve certain RHY demographics in NYC. 
We also defined an \textit{overflow shelter} that provides services to any RHY who could be placed in an existing RHY organization (that is, there is not a demographic mismatch), but the capacity within RHY organizations is insufficient to serve this youth.
\vspace*{-4mm}
\subsection{\emph{Runaway and Homeless Youth Profiles}} \label{s:youth}
To represent the demographics, needs, and desires of RHY in NYC we generated a simulated dataset of youth profiles using  primary and secondary data. Sources include data from NYC RHY programs, domain experts, and reports, which we used to generate proportions of varied youth demographics and needs in NYC.  These demographics are available in Table 4 and Table 5 of the Appendix. Once the proportions were established, we assigned different youth demographics, needs and desires to create a simulated dataset of RHY profiles to be used. Granular data on demographics, needs and desires of RHY are not publicly available and therefore, in our simulated dataset we assume all features in a youths profile are conditionally independent. 
\vspace*{-4mm}
\section{Modeling Capacity Expansion for RHY Organizations}\label{s:model}
We formulate the operational challenge of matching homeless youth to RHY organizations that provide housing and support services. Using mathematical optimization, we project the cost-minimizing capacity to meet the collective needs of youth.
\vspace*{-4mm}
\subsection{\emph{General Model Framework}} \label{s:framework}
Let {\boldmath{$T$}} be the set of days and {\boldmath{$S$}} the set of RHY organizations. We consider two types of RHY organizations serving homeless youth in NYC: (i) RHY shelters that offer housing (bed, food and basic necessities) and a variety of services in-house, denoted by {\boldmath{$S^b$}} $\subseteq$ {\boldmath{$S$}} and (ii) service providers that do not offer housing (beds) but provide other support services like medical assistance, mental health support and legal assistance (such as hospitals and mental health clinics) denoted by {\boldmath{$S^o$}} $\subseteq$ {\boldmath{$S$}}, where {\boldmath{$S^b$}} $\cup$ {\boldmath{$S^o$}} $=$ {\boldmath{$S$}} and {\boldmath{$S^b$}} $\cap$ {\boldmath{$S^o$}} $=$ $\emptyset$.  \\
 \indent Each RHY organization $s\in$ {\boldmath{$S^b$}} provides various services in-house, such as housing, healthcare and mental health support to youth with varying intensities. Let {\boldmath{$P$}} be the set of services provided to youth and  {\boldmath{$E_p$}} be the set of intensity levels available for each service $p\in$ {\boldmath{$P$}}. For most services, we consider 3 different intensities shown as $\epsilon \in$ {\boldmath{$E_p$}}: low, medium and high intensity depending on the amount of resources required to provide a particular service to youth. However, certain services have fewer intensity levels, for example, the housing service provided through RHY organizations has a single level of resource intensity, whereas mental health support services have 3 levels of intensity (low intensity: group therapy, medium intensity: weekly individual therapy, high intensity: seeing a psychiatrist and receiving medication). For ease of notation, we use index $i$ to represent every unique service-intensity pair $(p,\epsilon)$, where $i \in$ {\boldmath{$I$}} $=$ $\{(p,\epsilon) : p \in$ {\boldmath{$P$}} , $\epsilon \in$ {\boldmath{$E_p$}} $\}$. For brevity, we refer to each $i\in$ {\boldmath{$I$}} as a service, rather than its more precise ``service-intensity" name. A list of all services with their corresponding intensity options can be found in Table 2 of the Appendix. Figure \ref{fig:Organization} illustrates the two types of RHY organizations providing different services to RHY in NYC.
 
 \begin{figure}[htbp]
    \centering
    \includegraphics[width=1\textwidth]{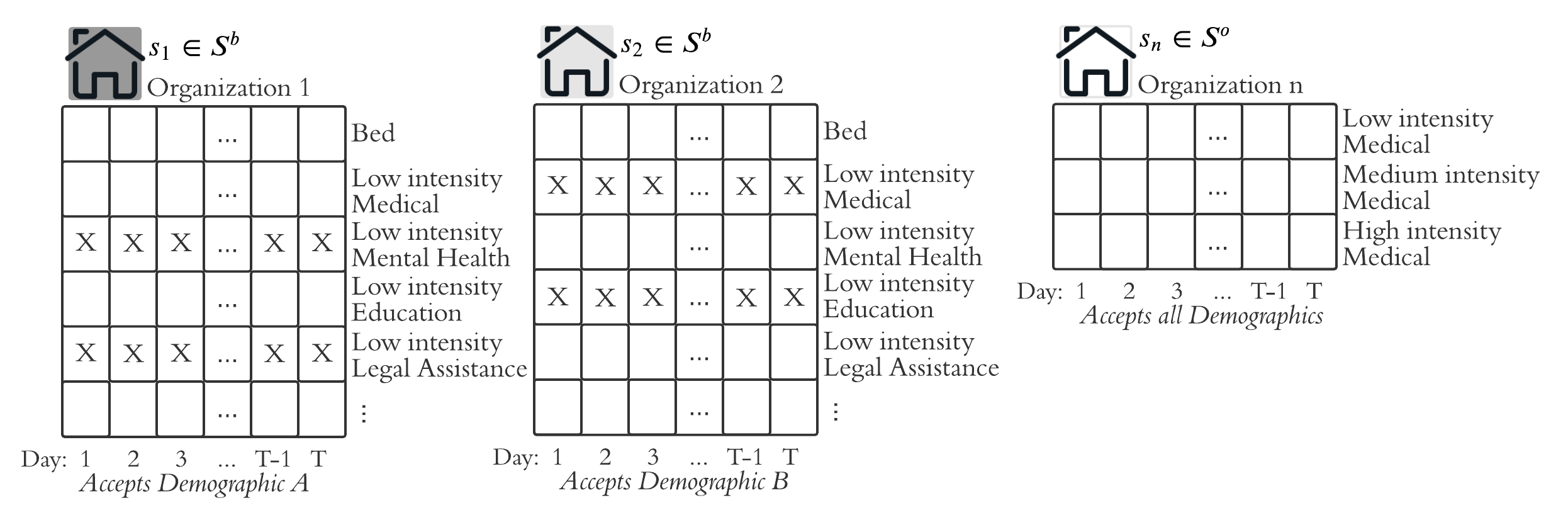}
    \caption{Illustrative example of sample RHY organizations that provide housing and support services to homeless youth, where X indicates the unavailable resources.}
    \label{fig:Organization}
\end{figure} 
 
\begin{figure}[h]
    \centering
    \includegraphics[width=0.7\textwidth]{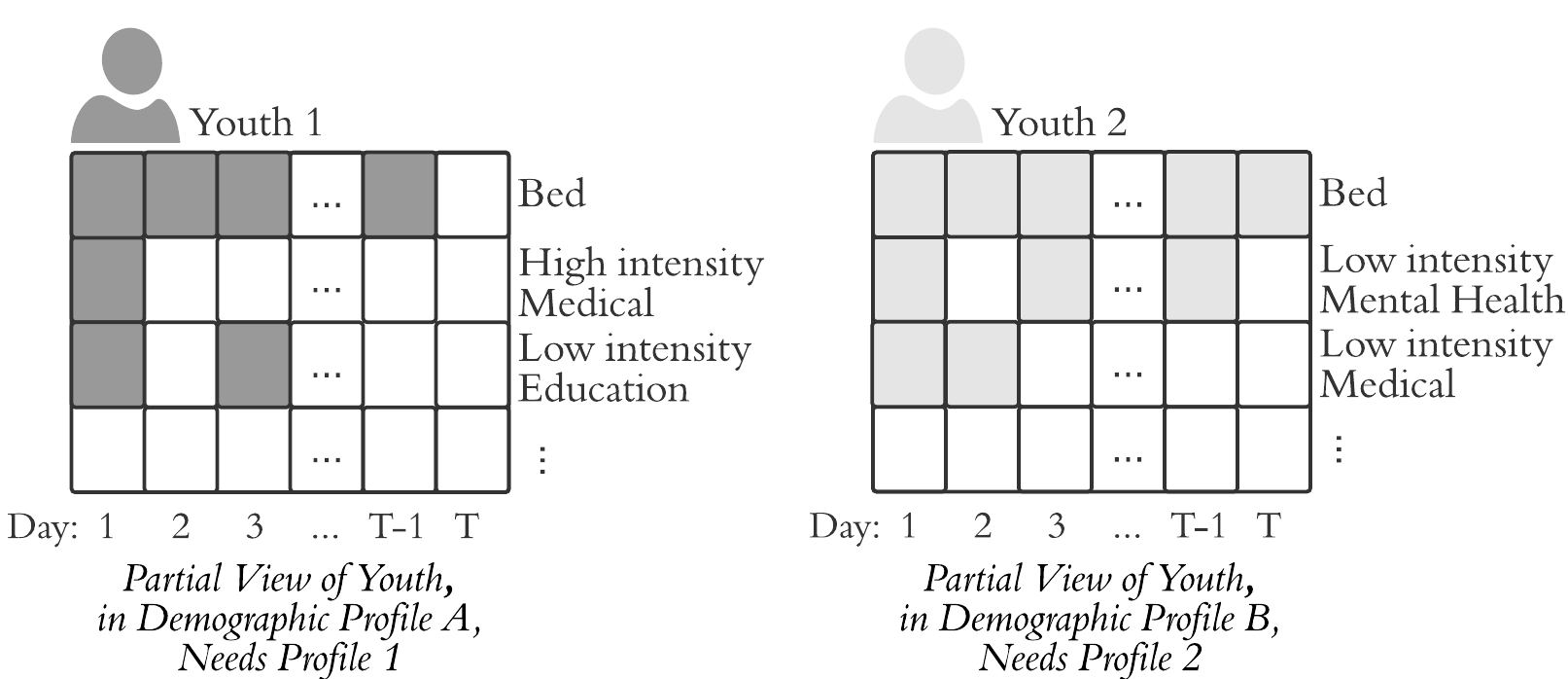}
    \caption{Examples of needs profiles for two youth who belong to two different demographics where shading indicates different demographics groups.}
    \label{fig:Youth}
\end{figure}

\begin{figure}[h]
    \centering
    \includegraphics[width=1\textwidth]{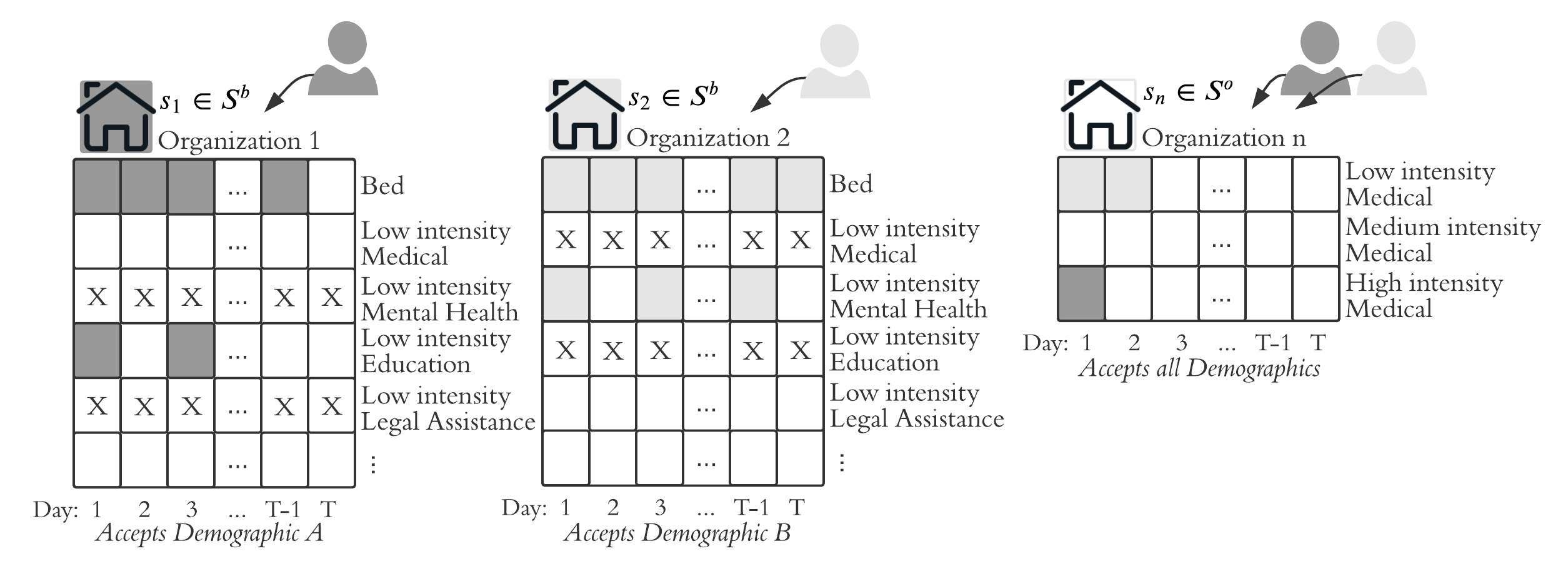}
    \caption{Youths are matched with RHY organization 1 and 2 respectively, considering the youth demographics and the accepted demographic at the RHY organizations. The services that are not provided in RHY organizations 1 or 2 are provided through $s_n \in$ {\boldmath{$S^o $}}. }
    \label{fig:Matching}
\end{figure}

Each RHY organization $s \in$ {\boldmath{$S^b$}} provides a variety of services to a subset of RHY demographics. For example, some organizations only provide services to females, some are specifically welcoming of LGBTQ+ people, some only serve youth under 21 years, and some do not accept families (RHY that have children of their own). These restrictions cause some youth to have particularly challenging time accessing services. To identify which demographics have reduced access, as discussed earlier in Section \ref{s:Provider}, we include an \textit{incompatibility set} as an additional RHY organization where youth who are unable to be placed in an existing RHY organization are assigned. \\
\indent Let {\boldmath{$Y$}} be the set of homeless youth in the system. Youth $y \in$ {\boldmath{$Y$}} arrives independently to the system on day $l_y \in$ {\boldmath{$T$}} with distinct need profile $\eta_y$ and demographics profile $\alpha_y$. Needs profile $\eta_y$ represents the needs that youth $y$ seeks from an RHY organization such as bed, financial assistance, or medical assistance, as well as the intensities and frequencies of these services. Each service $i$ in needs profile $\eta_y$ has a corresponding duration $d_{y,i}$ and frequency $f_{y,i}$. Demographic profiles $\alpha_y$ carry the age, gender, sexual orientation, child status, and HT victim information of youth $y$ and each attribute in the demographics profile is denoted as $n \in$ {\boldmath{$N$}}. For example $\alpha_y[1]$ is a binary value that denotes whether youth $y$ is a 16-year-old, and $\alpha_y[|${\boldmath{$N$}}$|]$ is a binary value that denotes whether youth $y$ has been a HT victim. An illustration of two distinct youth-needs profiles are depicted in Figure \ref{fig:Youth}. We use these needs profiles to match homeless youth with RHY organizations that serve their demographic and have available capacity to fulfill their unique needs. A representation of the matching process is depicted in Figure \ref{fig:Matching}. 

Sometimes RHY organizations have insufficient capacity at their facility to provide services to youth and must use creative approaches to meet demand such as by providing hotel vouchers. We model this capability through the concept of an \textit{overflow shelter}, which enables youth referral outside of the shelter. The overflow shelter captures the overall additional capacity needed in NYC after the in-house capacity of RHY organizations are extended as far as possible. Therefore, after the matching process, if the capacity within the RHY organization is insufficient to fulfill the needs of each youth, we project the additional capacity required by each in-house service, as well as how many youth should be directed to the overflow shelter. The reason we allow the capacity expansion through the overflow shelter is because RHY organizations have limited capacity within their facilities for additional resources, a phenomenon observed in areas where real estate is at a premium, such as NYC. Such an approach acknowledges the challenge of expanding in-house capacity by building another facility or expanding into a neighboring facility.
\vspace*{-4mm}
\subsection{\emph{Base Model: Capacity Expansion Optimization}} \label{s:basemodel}
We expand upon our explanation of the generalized model framework outlined in Section \ref{s:framework}. The sets, parameters and decision variables used in our optimization model are summarized in Tables \ref{t:sets}, \ref{t:parameters}  and \ref{t:decvariables}, respectively, and are further detailed in what follows. \\

\begin{table}[h!]
\begin{center}
\caption{List of sets used in mathematical modeling.}
\renewcommand{\arraystretch}{1.1}
\begin{tabular}{ p{1cm} p{13cm}  } 
 \toprule
 Symbol & Definition\\ [0.5ex] 
 \midrule
{\boldmath{$Y$}}                   & Set of youth in the system, indexed by $y$\\
{\boldmath{$S$}}                   & Set of RHY organizations in the system, indexed by $s$\\ 
{\boldmath{$T$}}                   & Time horizon over which services may be scheduled, indexed by $t$ \\
{\boldmath{$I$}}                   & Set of service-intensity pairs in system, indexed by $i$\\ 
{\boldmath{$N$}}                   & Set of demographics attributes in youth and organization profiles (age, gender, child and citizenship status, HT victim or survivor status), indexed by $n$\\
 \bottomrule
\end{tabular}
\label{t:sets}
\end{center}
\end{table}

\begin{table}[h!]
\begin{center}
\caption{List of parameters used in the mathematical modeling.}
\renewcommand{\arraystretch}{1.1}
\begin{tabular}{ p{1.5cm} p{12.5cm}  } 
 \toprule
 Symbol & Definition\\ [0.5ex] 
 \midrule
$l_y$                 & Arrival time of youth $y \in$ {\boldmath{$Y$}} to system\\
$d_{y,i}$             & Duration of service $i \in I$ for youth $y \in$ {\boldmath{$Y$}}\\
$f_{y,i}$             & Number of times service $i \in$ {\boldmath{$ I$}} should be provided to youth $y \in$ {\boldmath{$ Y$}}\\
$\omega_{y,i}$        & Time between appointments while providing service $i \in$ {\boldmath{$ I$}} to youth $y \in$ {\boldmath{$Y$}}\\
$[a_{y,i},b_{y,i} ]$  & Earliest and latest possible start times of service $i \in$ {\boldmath{$ I$}} to youth $y \in$ {\boldmath{$ Y$}} \\
$k_i$                 & Periodicity flexibility of service $i \in $ {\boldmath{$I$}} \\
$\alpha_y$            & Demographics profile of youth $y \in$ {\boldmath{$ Y$}}\\
$\beta_s$             & Demographics profile of RHY organization $s \in$ {\boldmath{$ S$}} \\
$\eta_y$             & List of requested services by youth $y \in$ {\boldmath{$ Y$}}\\
$\sigma_s$          & List of services offered by RHY organization $s \in$ {\boldmath{$ S$}}\\
$c_{s,i,t}$             & Capacity of service $i \in I$ in RHY organization $s \in$ {\boldmath{$ S$}} at time $t \in$ {\boldmath{$T$}}\\
$\mu_{s,i}$           & Maximum number of in-house resources that RHY organization $s \in$ {\boldmath{$ S$}} can have for service $i \in $ {\boldmath{$I$}} within their facility\\
$r_{y,s,i}$           & Cost of assigning youth $y \in$ {\boldmath{$ Y$}} to RHY organization $s \in$ {\boldmath{$ S$}} for service $i \in$ {\boldmath{$ I$}}  \\
$\gamma_{s,i}$        & Cost of allocating one unit of extra resource to RHY organization $s \in$ {\boldmath{$ S$}}, for service $i \in$ {\boldmath{$ I$}} \\
$\lambda_{s,i}$       & Cost of directing youth to overflow shelter from RHY organization $s \in$ {\boldmath{$ S$}} to receive service $i \in$ {\boldmath{$ I$}}\\
 \bottomrule
\end{tabular}
\label{t:parameters}
\end{center}
\end{table}

\begin{table}[h!]
\begin{center}
\caption{List of decision variables used in mathematical modeling.}
\renewcommand{\arraystretch}{1.1}
\begin{tabular}{ p{1.3cm} p{1.9cm} p{11cm} } 
 \toprule
Variable  & Type & Definition \\ [0.5ex] 
 \midrule
$U_{y,s,i}$         & Binary & Youth to RHY organization assignment decision variables are equal to 1 if youth $y \in$ {\boldmath{$  Y$}} is assigned to organization $s \in$ {\boldmath{$ S$}} to receive service $i \in$ {\boldmath{$  I$}}, 0 otherwise\\
$X_{y,s,i}^t$      & Binary & Time-dependent youth to RHY organization assignment decision variables are equal to 1 if youth $y \in$ {\boldmath{$ Y$}} is assigned to RHY organization $s \in$ {\boldmath{$ S$}} to receive service $i \in$ {\boldmath{$I$}} at time $t \in$ {\boldmath{$T$}}, 0 otherwise\\
$E_{s,i}^t$ & Continuous & Extra resource decision variables are equal to the amount of extra resources required to fulfill youths' collective demand at RHY organization $s \in$ {\boldmath{$ S$}}, to provide service $i \in$ {\boldmath{$ I$}} at time $t \in$ {\boldmath{$ T$}} \\
$O_{s,i}^t$ & Continuous & Overflow shelter decision variables are equal to the number of youths that are directed to the overflow shelter through RHY organization $s \in$ {\boldmath{$ S$}}, to receive service $i \in$ {\boldmath{$ I$}} at time $t \in$ {\boldmath{$ T$}}\\
 \bottomrule
\end{tabular}
\label{t:decvariables}
\end{center}
\end{table}

We model three important costs: (i) the cost of assigning youth to RHY organizations; (ii) the cost of adding extra resources to in-house services, and  (iii) the cost of directing youth to the overflow shelter. Assigning youth to RHY organizations that offer services in-house ($S^b$) has zero cost since we focus on the capacity expansion, and assignment to organizations that only offer support services such as hospitals and mental health clinics incur positive costs. These  costs are estimated using publicly available data sources such as hourly salary estimates and hotel voucher prices in NYC \citep{Gross-2021, NYCHRA-2022, Ziprecruiter-2022} and are sufficient to inform the capacity expansion required by each service provider to fulfill youth's needs \citep{Elluru-2019}. 

Our model is given by: \eqref{eq:obj}--\eqref{eq:periodicity}. Objective function \eqref{eq:obj} minimizes in aggregate the projected cost of: (i) assigning youth to RHY organizations, (ii) adding extra resources to RHY organizations for the services they provide, and (iii) directing youth to the overflow shelter.
\vspace{-5mm}
\begin{equation}
 Minimize \quad
 \sum_{y \in Y} \sum_{s \in S} \sum_{i \in I} \sum_{t \in  T}r_{y,s,i}  X_{y,s,i}^t +  \sum_{s \in S} \sum_{i \in I} \sum_{t \in T} \gamma_{s,i}  E_{s,i}^t + \sum_{s \in S} \sum_{i \in I} \sum_{t \in T} \lambda_{s,i} O_{s,i}^t \label{eq:obj}
\end{equation}
\vspace{-10mm}
\begin{subequations}
\begin{align}
& Subject \ to \  & &  \sum_{y \in Y} X_{y,s,i}^t \leq c_{s,i,t} +E_{s,i}^t +O_{s,i}^t,   &\forall s \in S, i \in I, t \in T,  \label{eq:capacity1} \\    
& & & c_{s,i,t} + E_{s,i}^t  \leq \mu_{s,i}, \quad & \forall s \in S, i \in I, t \in T, \label{eq:capacity2} \\
& & & \sum_{s \in S}U_{y,s,i} \leq 1, & \forall y \in Y, i \in I, \label{eq:capacity3}\\
& & & \sum_{t \in T}X_{y,s,i}^t \leq T \cdot U_{y,s,i}, &  \forall y \in Y, s \in S, i \in I, \  \label{eq:capacity4} \\
& & & \sum_{t \in T}\sum_{i \in I}X_{y,s,i}^t=0,   & \forall \{n\in N|  \eta_y[n] =1, \sigma_s[n]=0\},  y \in Y , s \in S, \label{eq:capacity5}\\
& & & U_{y,s,i} \in \{0,1\}, & \forall y \in Y , s \in S, i \in I, \label{eq:capacity6}\\
& & & X_{y,s,i}^t \in \{0,1\}, & \forall y \in Y , s \in S, i \in I, t \in T, \label{eq:capacity7} \\
& & & E_{s,i}^t \geq 0, \ O_{s,i}^t \geq 0, & \forall s \in S, i \in I, t \in T. \label{eq:capacity8and9} 
\end{align}
\label{eq:capacity}
\end{subequations}
\vspace{-10mm}

\indent The first system of constraints \eqref{eq:capacity} are related to the assignment of youth while considering the existing capacities within RHY organizations. The first constraint set \eqref{eq:capacity1} ensures that the number of youth assigned to a RHY organization to receive a service $i$ at shelter $s$ at time $t$ does not exceed the existing capacity of the service ($c_{s,i,t}$), extra in-house resources added to that service ($E_{s,i}^t$) and the number of youth directed to the overflow shelter ($O_{s,i}^t$). Constraint set \eqref{eq:capacity2} imposes an upper bound $\mu_{s,i}$ on the number of extra resources that can be added to an in-house service. Constraint set \eqref{eq:capacity3} ensures that a youth $y$ is receiving service $i$ from a single RHY organization at time $t$. Constraint set \eqref{eq:capacity4} ensures continuity of care; that is, it ensures that youth $y$ is receiving a service from a single RHY organization throughout the duration of their stay in the system. Finally, constraint set \eqref{eq:capacity5}  ensures that youth $y$ is not matched with RHY organization $s$, if RHY organization $s$ does not serve the demographic of youth $y$. Variable domains are stated in
 \eqref{eq:capacity6}-\eqref{eq:capacity8and9}.
\noindent \textbf{Service Delivery Time Windows} \\ Interviews with RHY organizations revealed that the timing of service provision (e.g. immediately upon arrival, in a few weeks) is nearly as important as the list of services provided to youth. All RHY organizations provide essential services that every youth must receive within the first 72 hours of arrival, such as: sexually transmitted disease testing, case management, and mental health assessment. We define the service delivery start time windows for each youth and service pair by $[a_{y,i}, b_{y,i}]$; the distributions they follow can be seen in Table 6 of the Appendix, and the system of constraints corresponding to service delivery time windows can be seen in \eqref{eq:time-window}. Constraint set \eqref{eq:tw1} requires a youth to be assigned a service between their earliest start time and latest end time. Constraint set \eqref{eq:tw2} imposes that the onset of provision of service $i$ to youth $y$ occurs between the earliest start time $a_{y,i}$ and latest start time $b_{y,i}$. 
\vspace{-2mm}
\begin{subequations}
\begin{align}
& & & \sum_{s \in S} \sum_{t = 0}^{a_{y_i} - 1} X_{y,s,i}^t + \sum_{s \in S} \sum_{t = b_{y,i} + d_{y,i} + 1}^{|T|} X_{y,s,i}^t = 0, & \ \forall y \in Y, i \in I, \label{eq:tw1}\\
& & & \sum_{s\in S}\sum_{t=a_{y,i}}^{b_{y,i}} X_{y,s,i}^t \geq 1, & \forall y \in Y , i \in I. \label{eq:tw2}
\end{align}
\label{eq:time-window}
\end{subequations}

\vspace{-8mm}
\noindent \textbf{Periodicity} \\
Real-world services may feature periodicity. Thus, the following constraint system \eqref{eq:periodicity} changes depending on whether the service is provided to youth in a periodic fashion ($i \in$ {\boldmath{$ I^{\omega}$}}) or not ($i \in$ {\boldmath{$ I^{n\omega}$}}). If the service is not periodic, constraint set \eqref{eq:period1} ensures that the summation over time of the time-dependent assignment variable $X_{y,s,i}^t$ equals the number of times youth requested the service $f_{y,i}$. In this case the time between the appointments is not significant. On the other hand, if the service is provided to youth in a periodic fashion, the time between the appointments should equal the periodicity $\omega_{y,i}$:
\begin{equation*}
\omega_{y,i}= \Bigl\lfloor  \frac{d_{y,i}}{f_{y,i}} \Bigr \rfloor. 
\label{eq:omega}
\end{equation*}

Thus, the periodicity $\omega_{y,i}$ is equal to the duration of service $d_{y,i}$ divided by the number of times the service is requested, $f_{y,i}$. To reflect operational reality, we also introduce parameter $k_i$ so as to allow flexibility in assigning youth to appointments within a small window of tolerance around $\omega_{y,i}$. For example, if service $i$ is provided to youth $y$ every Monday, with flexibility $k_i$ we ensure that youth $y$ can schedule their next appointment within $[-k_i, +k_i]$ days of Monday. Thus for a periodic service, constraint set \eqref{eq:period2} assures that summing $X_{y,s,i}^t$ over time equals the frequency $f_{y,i}$ while considering periodicity $\omega_{y,i}$ and flexibility $k_i$. Constraint set \eqref{eq:period3} imposes that youth $y$ is assigned to receive service $i$ at most once within the same periodicity flexibility window $[-k_i, +k_i]$.

\vspace{-8mm}
\begin{subequations}
\begin{align}
& & & \sum_{t=0}^{b_{y,i}+d_{y,i}} \sum_{s \in S}X_{y,s,i}^t=f_{y,i}, & \forall y \in Y, i \in I^{n\omega}, \label{eq:period1} \\
& & & \sum_{t=1}^{f_{y,i}} \sum_{s \in S} \sum_{k=-k_i}^{k_i} X_{y,s,i}^{t \cdot \omega_{y,i} +k}=f_{y,i}, & \forall y \in Y , i \in I^{\omega}, \label{eq:period2} \\
& & & \sum_{k=-k_i}^{k_i}X_{y,s,i}^{t \cdot \omega_{y,i} +k} \leq 1, & \forall y \in Y , \ s \in S \   i \in I^{\omega} , t \in \{0, ... , f_{y,i}\}. \label{eq:period3} 
\end{align}
\label{eq:periodicity}
\end{subequations}
\vspace*{-12mm}
\section{Computational Experiments and Results} \label{s:experiments}
We now present the computational setup for our experiments, the results obtained from  solving the mathematical models with varying model parameters, and discussions regarding  solution insights. All experiments were conducted using Gurobi Optimizer version 9.1 (2021) and Python 3.8.12, with up to 64 GB memory on an HPC cluster. Each instance was run with the Gurobi MIPGap optimality tolerance parameter set to 0.01. On average, it takes 24 hours to run a single scenario, which is reasonable as our model is used for long-term planning decisions. We conduct a variety of experiments on our simulated datasets using formulation \eqref{eq:obj}-\eqref{eq:periodicity}. 
\vspace*{-4mm}
\subsection{\emph{Computational Setup and Datasets}}
We evaluate the optimal capacity expansion for $|${\boldmath{$S$}}$| = 8$ RHY organizations that provide various TIL services at $|${\boldmath{$I$}}$| = 40$ intensity levels over $|${\boldmath{$T$}}$| = 180$ days (6 months). The increase of capacity at these organizations reduces a youth's vulnerability of being trafficked. During fiscal year 2021, 814 youth used approximately 300 TIL support program beds in NYC (986 and 1,221 youth were served in 2019 and 2020, respectively) \citep{NYC-gov}. The 8 RHY organizations that we collected primary and secondary data on provides 270 beds to youth. We thus assume in our base expansion model that on average $|${\boldmath{$Y$}}$| = 500$ youth arrive independently to these 8 RHY organizations to receive housing and (up to 40 different service-intensity pairs of) support services within 6 months. Interviews conducted with service providers revealed that approximately 90\% of existing resources are in use by RHY on any given day. Thus, while evaluating the capacity expansion, we assume that approximately 10\% of the current resources are idle for any of the incoming 500 youth to use.\\
\indent In light of the uncertainty regarding RHY data, we performed extensive sensitivity analyses around several key model parameters to determine their effect on the optimal capacity expansion for RHY organizations. We changed one parameter at a time in the base model to test the impact of four additional levels of that parameter. We also considered an additional COVID-19 scenario where two parameters change simultaneously. Table \ref{t:experiments} summarizes the ranges of values used in our sensitivity analysis, in total covering 14 scenarios (1 base case; 12 sensitivity analyses; and 1 COVID-19 scenario). We ran each  scenario 10 times to consider the 10 unique youth-needs profiles, yielding 140 runs. We then present our results and insights of sensitivity analyses in Section \ref{s:Sensitivity}. 
\begin{table}[]
\caption{Experimental parameters varied; bolded values indicate base expansion model parameters.}
\centering
\begin{tabular}{lcc}
\toprule
Parameter                       & Symbol      & Level                                 \\ 
\midrule
Number of youth                 & {\boldmath{$Y$}} & 400, 450, \textbf{500}, 550, 600                          \\
Abandonment percentage of youth & $\theta$ & 10\%, 15\%, \textbf{20\%}, 25\%, 30\%                      \\
Duration of service             & $d_{y,i}=N(60,15)$   & 0.8$d_{y,i}$, 0.9$d_{y,i}$, {\boldmath{$d_{y,i}$}}, 1.1$d_{y,i}$, 1.2$d_{y,i}$ \\
COVID-19 effect                & {\boldmath{$Y$}}$, \ c_{s,i,t}$  & $|${\boldmath{$Y$}}$| = 400$ $, \ c_{s,i,t}'=0.5c_{s,i,t}$ \\
\bottomrule
\end{tabular}
\label{t:experiments}
\end{table}

\subsection{\emph{Baseline Model Results and Insights}}
This section discusses our base expansion model results, which provide the optimal capacity expansion required by RHY organizations considering our base input parameters. The capacity expansion we present assumes that youth are matched with RHY organizations as efficiently as possible, thus, the capacity expansion we present is likely to represent a conservative estimate on the extra resources required by RHY organizations to fulfill collective youth needs. To capture the effect of variation, we created 10 different youth-needs profiles ($\eta_y$) using the information given in Table 5.
\vspace*{-4mm}
\subsubsection{\emph{Bed expansion}}
We allow for four different bed types: (i) existing in-house beds, (ii) extra in-house beds that expand the capacity within the facility, (iii) overflow shelter beds when bed space in the facility has been maximally expanded, and (iv) incompatibility set ({\boldmath{$\Psi$}}) beds when none of the existing 8 RHY organizations are able to serve a particular youth's demographic. The demographics each RHY organization serves, services they provide and the number of existing beds appear in Tables 1 and 3 of the Appendix. Figure \ref{fig:Beds_organizations} depicts the average number of youth across 10 scenarios receiving these different bed types from 8 RHY organizations in the optimal solution to our base model.\\
 \begin{figure}[htbp]
    \centering
    \includegraphics[width=0.75\textwidth]{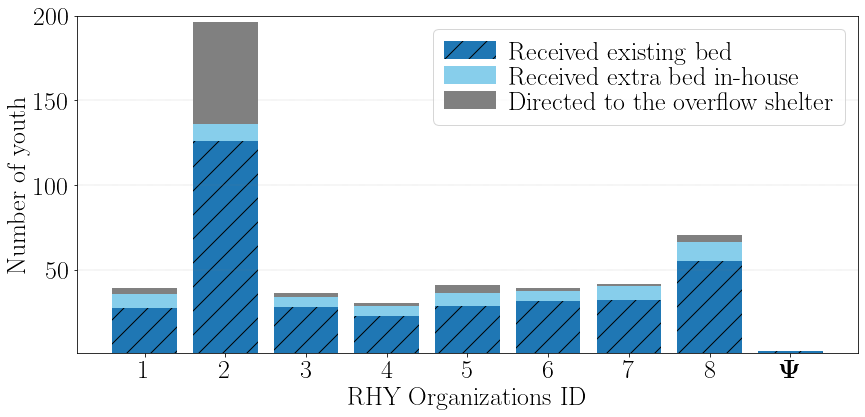}
    \caption{Average number of beds received by type, across 10 runs, 8 organizations, and 6 months.}
    \label{fig:Beds_organizations}
\end{figure}

\begin{figure}[htbp]
    \centering
    \includegraphics[width=0.8\textwidth]{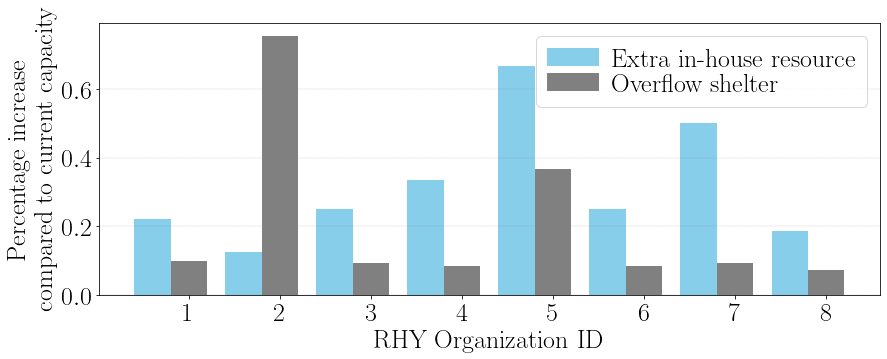}
    \caption{Per RHY organization, the average percentage increase in extra in-house resources and number of youth directed to overflow shelter compared to existing in-house capacity over 6 months. }
    \label{fig:Beds_percentages}
\end{figure}

\indent All of the RHY organizations in the system used their existing capacity and added the maximum possible amount of extra beds, yet still needed to resort to creative options for beds, accommodating youth via the overflow shelter. Such results underscore the critical need for more RHY organizations in NYC. Out of 500 youth, on average across 10 runs, 63 youth received an extra in-house bed, 80 youth had to receive a bed through the overflow shelter due to capacity restrictions within existing facilities, and 2 youth had to be placed in the incompatibility set ({\boldmath{$\Psi$}}) as none of the 8 RHY organizations could serve their particular demographic. The model demonstrates that simply adding more capacity to  existing shelters is insufficient as certain demographics experience access challenges. Figure \ref{fig:Beds_percentages} represents the percentage of capacity expansion required by each RHY organization. The average total percentage increase in resources needed for extra in-house beds and overflow shelters, collectively, across the eight organizations is approximately 52\%.

Figure \ref{fig:Beds_organizations} shows that Organization 2 requires the greatest amount of additional capacity. This  particular organization has fewer entry restrictions compared to other organizations, thus can serve a broader range of populations. Accordingly, the optimal solution increases the capacity of Organization 2 more than any other shelter. Notably, this organization already has 80 TIL beds within their facility. Thus, considering the youth who are directed to the incompatibility set {\boldmath{$\Psi$}} and  the significant capacity expansion needs in Organization 2, our results show the impact of collective entry requirements (such as age and gender restriction) on meeting the needs of NYC RHY. 
\vspace*{-4mm}
\subsubsection{\emph{Support Service Expansion}}
\begin{figure}[h]
    \centering
    \includegraphics[width=0.8\textwidth]{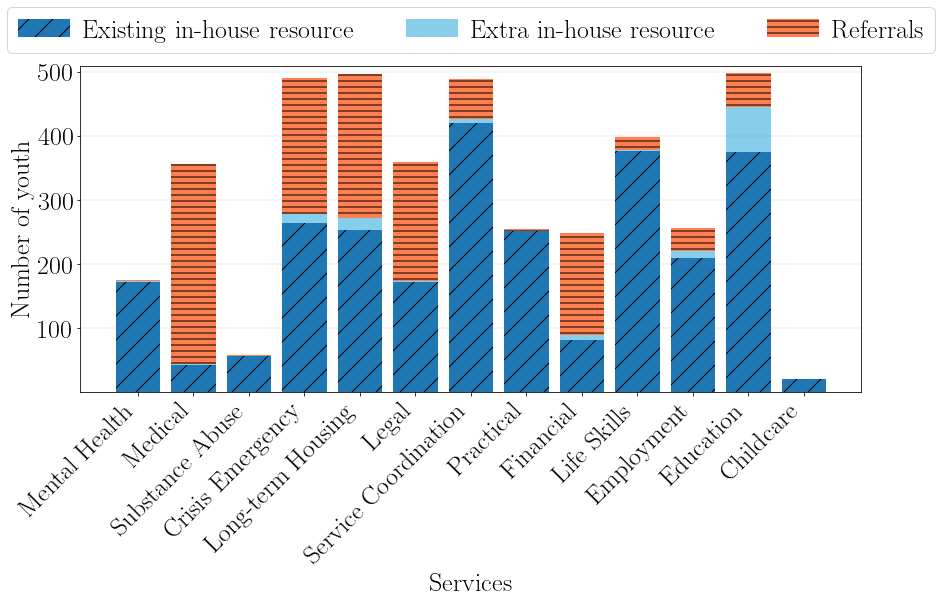}
    \caption{Average number of support services received, by type, over 6 months, across 10 runs.}
    \label{fig:Services_organizations}
\end{figure}

Recall that RHY organizations provide support services to youth (i) with existing in-house resources, (ii) with extra in-house resources added to expand the capacity within the facility, and (iii) through referrals ({\boldmath{$S^o$}}). Figure \ref{fig:Services_organizations} illustrates the average number of youth across 10 runs who received these services by source (note that the intensities of services are combined). In our simulated data, medical assistance is required by the majority of youth, which is currently provided primarily through referrals. Additionally, a large portion of crisis emergency, long-term housing, legal, and financial assistance must be provided through referrals, while educational assistance requires additional in-house resources as it is less costly to provide. When considering the addition of resources, referrals are preferred over in-house services as the latter is more costly. However, as Figure \ref{fig:Services_organizations} illustrates, there is a constant need for referrals for specific services: medical, crisis emergency, long-term housing, legal, and financial assistance. This underscores that adding extra in-house resources would increase convenience and access for youth, even if more costly. 
\begin{figure}[h]
    \centering
    \includegraphics[width=0.8\textwidth]{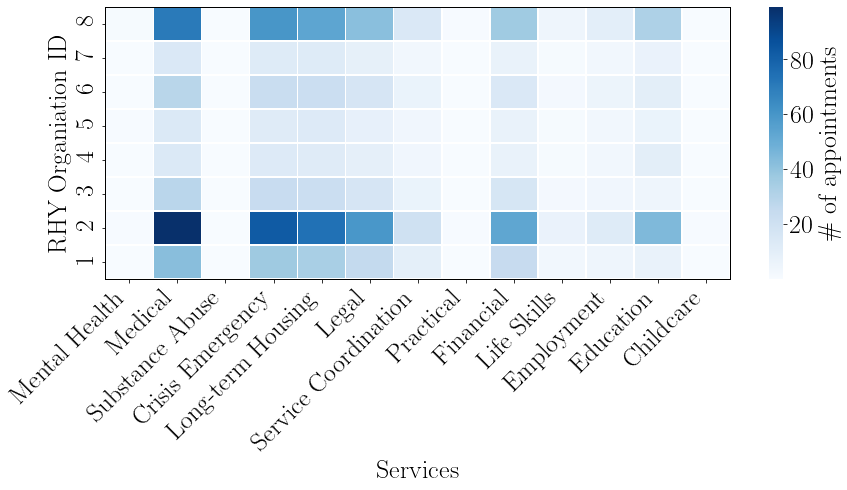}
    \caption{Average number of extra hours required, per RHY organization and support service, to fulfill youth demand, across 10 runs (note that extra in-house resources and referrals are combined). }
    \label{fig:Services_Heatmap}
\end{figure}

Figure \ref{fig:Services_Heatmap} depicts the number of extra appointments required by RHY organizations to provide support services. These results show that the additional resources needed for medical, legal, financial and educational assistance are significant in reducing vulnerability to trafficking. Figure \ref{fig:Services_Heatmap} shows that the number of youth requesting mental health support, substance abuse, childcare and practical assistance is lower than other services. This reflects perceived youth needs and mirrors insights gleaned from key stakeholders. RHY may hesitate requesting mental health and substance abuse support due to judgment, doubt, fear, and misinformation \citep{Donley-2021}.

\vspace*{-4mm}
\subsection{\emph{Sensitivity Analysis Insights}}\label{s:Sensitivity}
We now discuss sensitivity analysis insights gained by varying model parameters listed in Table \ref{t:experiments}, including changes in the optimal capacity expansion plan and RHY organization types to expand. 
\vspace*{-4mm}
\subsubsection{\emph{Arrival Rate}}
The number of RHY seeking TIL opportunities is highly correlated with extraneous factors such as weather, safety concerns, and the political environment \citep{Tambe-2018}. Accordingly, we perform a sensitivity analysis on youth arrival rate, using data from our service provider interviews.
\vspace{-6mm}
\begin{figure}[htbp]
    \centering
    \includegraphics[width=0.8\textwidth]{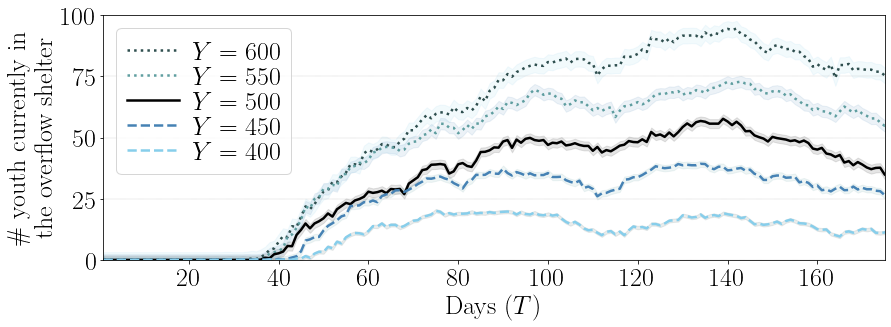}
    \caption{The effect of total number of youth in the system, on the daily trend in number of youth at the overflow shelter referred from Organization 2 over 6-months, with 90\% confidence intervals over 10 runs.}
    \label{fig:arrival}
\end{figure}

All eight RHY organizations show a similar capacity expansion behavior, which we illustrate with the largest shelter, Organization 2. Figure \ref{fig:arrival} shows the number of youth directed to receive an overflow shelter from Organization 2, and up to 10 beds can be added within the facility ($E^t_{2,1}=10$). In the scenario with the fewest youth ($|${\boldmath{$Y$}}$| = 400$), across 10 runs, on average 31 youth remain unable to receive services with  existing resources from Organization 2. Unsurprisingly, the number of youth being directed to the overflow shelter peaks when $|${\boldmath{$Y$}}$| = 600$; even when 10 more beds are added in-house, on average 95 youth are unable to receive a bed within the facility. 
\begin{table}[h] 
    \caption{System overflow metrics comparing number of youth arriving with the base expansion model of 500 youth considering the average across 10 runs.}
    \centering
    \begin{tabular}{p{4.3cm} p{0.4cm} p{0.4cm} p{0.4cm}p{0.4cm}p{0.4cm}}
\toprule
\multirow{2}{4.4 cm}{System Overflow Metrics}& \multicolumn{5}{c}{Number of Youth} \\ \cline{2-6}
& \multicolumn{1}{c}{{400}} & \multicolumn{1}{c}{{450}} & \multicolumn{1}{c}{\textbf{500} } &
\multicolumn{1}{c}{{550}} & \multicolumn{1}{c}{{600}} \\ \hline
{Max Overflow Beds} & \multicolumn{1}{c}{37.0} &  \multicolumn{1}{c}{52.0} &  \multicolumn{1}{c}{\textbf{87.0}} & \multicolumn{1}{c}{90.0}   & \multicolumn{1}{c}{110.0}\\
{Mean Overflow Beds} & \multicolumn{1}{c}{33.7				
} & \multicolumn{1}{c}{46.7} & \multicolumn{1}{c}{\textbf{80.4}} & \multicolumn{1}{c}{82.8}  & \multicolumn{1}{c}{102.0}\\
\midrule
{Overflow Cost Change} & \multicolumn{1}{c}{-65\%} & \multicolumn{1}{c}{-29\%} & \multicolumn{1}{c}{\textbf{0\%} } & \multicolumn{1}{c}{+35\%} &\multicolumn{1}{c}{+67\%} \\ 
{Referral Cost Change} & 		
\multicolumn{1}{c}{-20.1\% 		
} & 
\multicolumn{1}{c}{-39.3\%} & \multicolumn{1}{c}{\textbf{0\%}}  &
\multicolumn{1}{c}{+72.9\%} &  \multicolumn{1}{c}{+94.5\%}  \\
\bottomrule
\end{tabular}
    \label{tab:Arrival_overflow}
\end{table}

From a system-wide perspective, our base expansion model assumes that 500 youth arrive to the system over a six-month period. Considering alternative arrival numbers, in Table \ref{tab:Arrival_overflow} we see that the six-month average of the overflow shelter need increases by 22\% (from 80.4 to 102) when the number of youth increases by 20\% ($|${\boldmath{$Y$}}$|$ $= 600$); and decreases by 58\% (from from 80.4 to 33.7) when the number of youth decreases by 20\% ($|${\boldmath{$Y$}}$| = 400$). Quantifying the extent to which the arrival rate affects capacity is critical supporting data for capacity expansion decisions.
\vspace*{-4mm}
\subsubsection{\emph{Duration of Service}}
Service provider interviews revealed average lengths of stay for youth of nearly 60 days. Yet, as discussed in Section \ref{s:Provider}, the service duration $d_{y,i}$ for youth $y$ depends on the individual youth, organization and service type. Thus, we vary service duration $d_{y,i}$ which in our base expansion model follows a normal distribution with mean of 60 days and standard deviation of 15 days. 

\begin{figure}[htbp]
    \centering
    \includegraphics[width=0.8\textwidth]{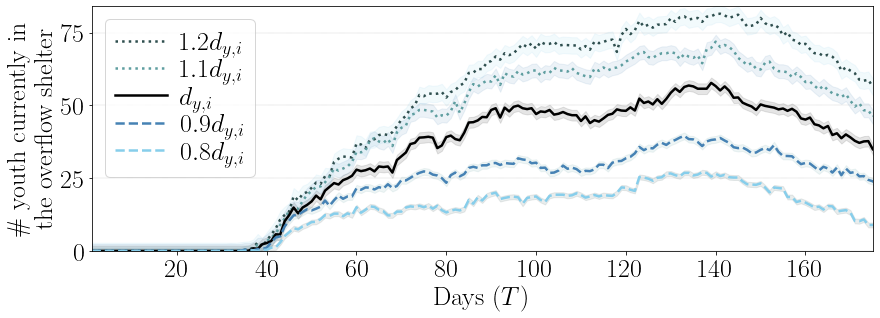}
    \caption{The effect of average length of stay, on the daily trend in number of youth at the overflow shelter referred from Organization 2 over 6 months, with 90\% confidence intervals over 10 runs.} 
    \label{fig:duration}
\end{figure}

\begin{table}[htbp] 
    \caption{Average system overflow metrics under base expansion model varying $d_{y,i}$ levels.}
    \centering
    \begin{tabular}{p{4.3cm} p{0.4cm} p{0.4cm} p{0.4cm}p{0.4cm}p{0.4cm}}
\toprule
\multirow{2}{4.4 cm}{System Overflow Metrics}& \multicolumn{5}{c}{Duration of Service} \\ \cline{2-6}
& \multicolumn{1}{c}{{$0.8d_{y,i}$}} & \multicolumn{1}{c}{{$0.9d_{y,i}$}} & \multicolumn{1}{c}{\boldmath{$d_{y,i}$} } &
\multicolumn{1}{c}{{$1.1d_{y,i}$}} & \multicolumn{1}{c}{{$1.2d_{y,i}$}} \\ \hline
{Max Overflow Beds} & \multicolumn{1}{c}{56.0} &  \multicolumn{1}{c}{62.0} &  \multicolumn{1}{c}{\textbf{87.0}} & \multicolumn{1}{c}{97.0}   & \multicolumn{1}{c}{117.0}\\
{Mean Overflow Beds} & \multicolumn{1}{c}{48.0} & \multicolumn{1}{c}{	58.7} & \multicolumn{1}{c}{\textbf{80.4}} & \multicolumn{1}{c}{96.0}  & \multicolumn{1}{c}{114.0}\\
\midrule
{Overflow Cost Change} & \multicolumn{1}{c}{-58\%} & \multicolumn{1}{c}{-31\%} & \multicolumn{1}{c}{\textbf{0\%} } &
\multicolumn{1}{c}{23\%} &\multicolumn{1}{c}{+50\%} \\ 
{Referral Cost Change} & 		
\multicolumn{1}{c}{-26\%} & 
\multicolumn{1}{c}{-12\%} & \multicolumn{1}{c}{\textbf{0\%}}  &
\multicolumn{1}{c}{+12\%} &  \multicolumn{1}{c}{+37\%}  \\
\bottomrule
\end{tabular}
    \label{tab:duration_overflow}
\end{table}
Variation in duration $d_{y,i}$ also affects the number of times a youth needs each service ($f_{y,i}$), thereby influencing the cost-minimizing capacity expansion. As seen in Figure \ref{fig:duration}, decreasing the average length of stay by 20\% for 500 youth eliminates 40\% (from 80.4 to 48) of the need for the overflow shelter. However, reduced length of stays at a shelter is likely to disrupt the rehabilitation process and increase future vulnerability and recidivism for those at risk of trafficking. On the other hand, in Table \ref{tab:duration_overflow} it is shown that a 20\% increase in the duration results in 117 youth on average across 10 runs, being unable to access an existing bed, raising their vulnerability to trafficking. 

A shorter service duration results in reduced overflow costs, slightly lower referral costs, and less capacity expansion as revealed in Table \ref{tab:duration_overflow}. There is a trade-off when considering 
reducing the service duration for youth; while this may alleviate some current capacity limitations and temporarily improve access to housing and support services, it will likely disrupt the much needed efforts for youth to be able to successfully exit trafficking and exploitative experiences. Capacity expansion plans should consider the average service duration of youth, and should provide for the possibility of extending the stay of youth until a safe and stable living arrangement is identified. 
\vspace*{-4mm}
\subsubsection{\emph{Abandonment Percentage}} 
A youth may stop receiving a particular service or abandon the system completely for various reasons, such as feeling limited by organizational restrictions, avoiding conflict and abuse, relapse, health concerns, or the inability to find stable housing \citep{Donley-2021}. Uncertainty exists regarding the number of youth abandoning the system while receiving services. Accordingly, we performed sensitivity analysis regarding the percentage of youth abandoning the system ($\theta$), which in our base model is set to 20\%. Service provider interviews revealed that almost half of youth who abandon the system leave within the first three days of arrival. Therefore when modeling abandonment, we assume that half of the youth who abandon have a service duration ($d_{y,i}$) that follows a normal distribution with mean of 3 days and standard deviation of 0.5 days, whereas the other half has their duration decrease to third of the original duration. 

\begin{figure}[htbp]
    \centering
    \includegraphics[width=0.8\textwidth]{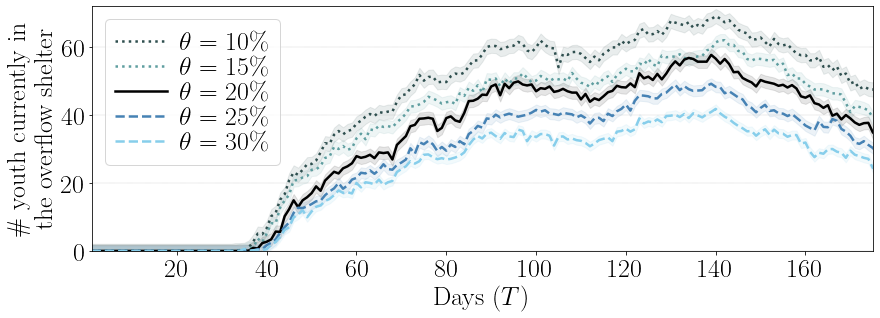}
    \caption{The effect of abandonment, on the daily trend in the number of youth at the overflow shelter referred from Organization 2 over 6 months, with 90\% confidence intervals over 10 runs.} 
    \label{fig:abandonment}
\end{figure}

Table \ref{tab:abandonment_overflow} shows that a 10\% increase in abandonment decreases the average need for overflow by 20\% (from 80.4 to 64.7). In comparison, a 10\% decrease in abandonment increases the average overflow by 16\% (from 80.4 to 93). When the percentage of youth abandonment is between 10--30\%,  Figure \ref{fig:abandonment} shows that in Organization 2 the in-house capacity and extra in-house resources are exhausted, resulting in the need for overflow around day 40. When the abandonment level is at 10\%,  the average number of referred youth actively in the overflow shelter peaks at 70. In contrast, on average across 10 runs, only 42 youth  are directed when abandonment increases to 30\%. 

\begin{table}[htbp] 
    \caption{Average system overflow metrics under base expansion model, varying abandonment levels.}
    \centering
    \begin{tabular}{p{4.3cm} p{0.4cm} p{0.4cm} p{0.4cm}p{0.4cm}p{0.4cm}}
\toprule
\multirow{2}{4.4 cm}{System Overflow Metrics}& \multicolumn{5}{c}{Abandonment \%} \\ \cline{2-6}
& \multicolumn{1}{c}{{30\%}} & \multicolumn{1}{c}{{25\%}} & \multicolumn{1}{c}{\textbf{20\%} } &
\multicolumn{1}{c}{{15\%}} & \multicolumn{1}{c}{{10\%}} \\ \hline
{Max Overflow Beds} & \multicolumn{1}{c}{67.0} &  \multicolumn{1}{c}{75.0} &  \multicolumn{1}{c}{\textbf{87.0	}} & \multicolumn{1}{c}{92.0}   & \multicolumn{1}{c}{97.0}\\
{Mean Overflow Beds} & \multicolumn{1}{c}{64.7} & \multicolumn{1}{c}{69.7} & \multicolumn{1}{c}{\textbf{80.4}} & \multicolumn{1}{c}{88.2}  & \multicolumn{1}{c}{93.0	}\\
\midrule
{Overflow Cost Change} & \multicolumn{1}{c}{-27\%} & \multicolumn{1}{c}{-15\%} & \multicolumn{1}{c}{\textbf{0\%} } &
\multicolumn{1}{c}{14\%} &\multicolumn{1}{c}{+29\%} \\ 
{Referral Cost Change} & 		
\multicolumn{1}{c}{-18\%} & 
\multicolumn{1}{c}{-8\%} & \multicolumn{1}{c}{\textbf{0\%}}  &
\multicolumn{1}{c}{+10\%} &  \multicolumn{1}{c}{+16\%}  \\
\bottomrule
\end{tabular}
    \label{tab:abandonment_overflow}
\end{table}

While it is challenging to estimate the rate and drivers of abandonment, these appear to have significant influence on capacity expansion decisions. Although beyond the scope of this study, we maintain that it is crucial to understand the reasons for youth abandonment and to decrease the percentage of youth abandoning the system. In turn, reduced abandonment is likely to decrease the vulnerability of runaway and homeless youth to being trafficked or exploited. 
\vspace*{-4mm}
\subsubsection{\emph{The Effect of COVID-19}}{\label{s:covid}} 
Access to scarce healthcare resources is a common challenge during the COVID-19 pandemic, especially for vulnerable populations \citep{Dorn-2020}. We evaluate the healthcare needs of RHY as access to healthcare would reduce their vulnerability to HT \citep{Duncan-2019}. Over the course of the COVID-19 pandemic, semi-regular meetings with NYC stakeholders revealed that the number of youth arriving to the RHY organizations decreased due to concerns of infection. RHY organizations also experienced significant staff unavailability and attrition due to positive COVID-19 cases and various vaccine mandates. To evaluate the effect of COVID-19 on the need for capacity expansion, we decrease the existing capacity of in-house services to 50\% of the base model and decrease the number of youth arriving to $Y=400$, in accordance with stakeholder recommendations and publicly available data \citep{NYC-gov}. 
\begin{figure}
  \centering
  \begin{subfigure}[t]{0.53\textwidth}
    \centering
    \includegraphics[width=1\textwidth]{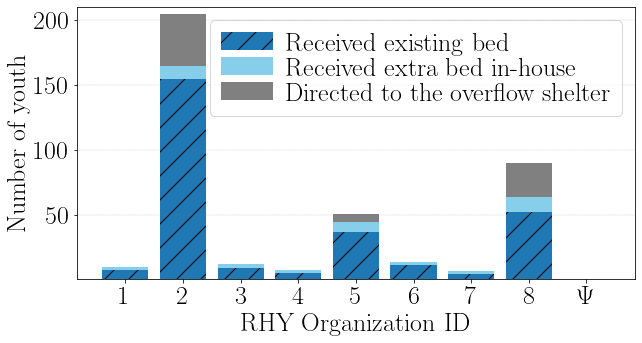}
    \caption{Beds received from varying sources.}
    \label{fig:Beds_covid}
  \end{subfigure}
  \begin{subfigure}[t]{0.46\textwidth}
    \centering
    \includegraphics[width=1\textwidth]{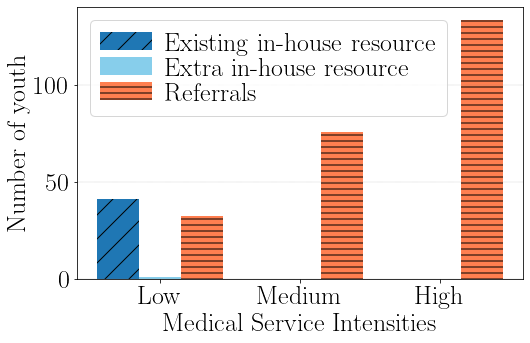}
    \caption{Medical services received from varying sources.}
    \label{fig:services_covid}
  \end{subfigure}  
  \caption{Average number of beds and different intensity medical resources received by youth in organizations considering COVID-19, across 10 runs and 6 months.}
  \label{fig:covid}
\end{figure}

While the halving of in-house capacity left more physical space for extra in-house resources, staffing beds during the pandemic remained a challenge. Figure \ref{fig:Beds_covid} shows the need to direct youth to overflow shelter beds still exists; averaged across 10 runs, 113 of the 400 youth did not receive an existing bed. Figure \ref{fig:services_covid} shows that in an optimal expansion, averaged over 10 runs, only 14\% of youth could be accommodated with existing in-house medical resources; and while an additional 1\% could be accommodated via extra in-house resources, a full 85\% must have their medical needs met via referrals. However, during the pandemic, referrals to larger healthcare institutions became less desirable due to health concerns \citep{Li-2020} and more difficult to access due to redirected healthcare resources and social distancing. Such a situation reveals the need for additional in-house medical resources. Our model can evaluate the trade-off between the increased costs of providing in-house services and access to these services in situations such as a pandemic.
\vspace*{-4mm}
\section{Conclusion}\label{s:conclusion}
Operations research and analytics research efforts to disrupt the supply of sex and labor trafficking victims are relatively new \citep{Dimas-2022}. We use optimization-based techniques to evaluate  cost-minimizing capacity expansion options for shelters serving  runaway and homeless youth and young
adults (RHY), which we believe is a first such attempt to address the vulnerability of a population extremely susceptible to trafficking \citep{Wright-2021, Hogan-2020}. We present an integer linear programming model that incorporates stochastic youth arrivals and length of stay; youth abandonment;  service delivery time windows, as well as periodic and non-periodic services.  
Through careful variable definition, we allow for three types of time-dependent capacity expansion: adding extra resources to in-house services, directing youth to an overflow shelter when there is no longer an ability to expand in-house capacity, and an \textit{incompatibility set} that serves youth who are otherwise unable to receive housing and support services due to demographic mismatch. While we illustrate our approach with a case study aimed at expanding transitional independent living service capacity for RHY organizations in New York City (NYC), the same (or a similar) model could be readily employed for non-profit shelter organizations in other locations. \\
\indent We discovered new insights that have the potential to impact human trafficking disruption efforts by increasing the accessibility of services aimed at reducing youth vulnerability of being trafficked and exploited. Our systematic and data-driven analysis addresses extremely resource-constrained contexts by providing a capacity expansion strategy for NYC with practical impact. Overall, our study represents an innovative use of mathematical modeling to address a capacity expansion problem with a broader societal impact.\\
\indent As our mathematical model uses optimization to most efficiently place youth in RHY organizations, it is worthwhile to consider that the expansion we recommend is likely to represent a conservative estimate on the extra resources required by RHY organizations to fulfill the collective needs of RHY. Moreover, through a thorough sensitivity analysis on key model parameters, we demonstrate that, even in the least demanding scenario, there is need to expand in-house capacity to its upper bound and add 20 overflow beds. Although our model includes capacity expansion to existing RHY organizations (as opposed to adding completely new organizations), when we look at the scenarios with higher demands, the results clearly indicate that more shelters are needed in NYC. Interestingly, the sensitivity analyses show that change in service duration affects the overall capacity expansion more than the change in arrival rate or the abandonment rate. We see that a 20\% increase in service duration, and arrival rate increases the need for overflow beds by on average 35\% and 27\% across 10 runs, respectively. In comparison, a 20\% decrease in abandonment rate, which results in more youth using services, only increases the need for overflow beds by 12\%. Additionally, we substantiate that inclusion and exclusion criteria dictating who may receive shelter services complicate access to housing resources for certain demographics. 
The aforementioned insights were informed through analysis of model output based on estimated RHY profiles using a mix of primary and secondary data. While a potential limitation, we believe our secondary sources have provided a reasonable representation of the current demographic and needs profiles of youth, and moreover a solid foundation for future work. Our results can be used as a template for future analysis, including the analysis of primary data collection on RHY demographics, needs and desires.\\ 
\indent While this study specifically focuses on capacity expansion of RHY organizations within NYC, there remains great potential for organizations in other locations provided sufficient data exists. Further extensions include optimally locating the additional overflow shelter and extending the model to allow RHY organizations to share resources with one another. Additionally, as demand for housing and support services greatly exceeds the existing capacity, the optimal deployment will require more capacity than is feasible to add at one time. Therefore, identifying an actionable capacity expansion plan that details how to implement the capacity deployment over time would offer additional utility. Moreover, there remains potential to embed our optimization approach into a decision-support tool to further facilitate the decision-making process in NYC.\\
\indent In conclusion, this study provides preliminary evidence of the value of incorporating the preferences and needs of vulnerable populations into humanitarian operations research problems. Our approach benefits government and nonprofit decision-makers by offering a means to effectively evaluate the allocation and expansion of scarce resources, while readily enabling sensitivity analyses to examine the effect of demand changes on the optimal expansion of housing and support services.
\vspace*{-4mm}
\section{Data Availability, and Acknowledgements}
The authors confirm that the aggregated data supporting the findings of this study are available within the article and its supplementary materials. As stated in the Institutional Review Board, due to the nature of this research, study participants (such as RHY organizations) did not agree for their data to be shared publicly; thus detailed, identifiable supporting data is not available.
This work is supported by the National Science Foundation under Grant No. CMMI-1935602. 
The authors would also like to thank the New York City Mayor's Office and New York Coalition for Homeless Youth for their insights and all of the RHY service organizations that helped inform this study. Special thanks to Andrea Hughes for her time and efforts. 

\bibliographystyle{chicago}
\spacingset{1}
\bibliography{IISE-Trans}

\clearpage
\section*{Appendix}
The RHY and RHY organization demographic, needs and service profiles are shown in Appendix. \\

\begin{table}[htbp]
\centering
\caption{RHY organization demographic profiles used in modeling ($\beta_s$), where 1 indicates the RHY organization accepts a particular demographic}
\begin{tabular}{clllllllll}
\toprule
\multicolumn{2}{c}{\multirow{2}{*}{Demographic Attributes}} & \multicolumn{8}{c}{Organization} \\ \cline{3-10} 
\multicolumn{2}{c}{}                                        & 1  & 2  & 3  & 4 & 5 & 6 & 7 & 8 \\ \midrule
\multirow{2}{*}{Age}                & 21 below              & 1  & 1  & 1  & 1 & 1 & 1 & 1 & 1 \\
                                    & 21+                   & 1  & 0  & 1  & 1 & 0 & 0 & 1 & 1 \\ \midrule
\multirow{7}{*}{Gender}             & Cis-gender Male       & 0  & 1  & 1  & 0 & 1 & 1 & 1 & 1 \\
                                    & Cis-gender Female     & 0  & 1  & 1  & 1 & 1 & 1 & 1 & 1 \\
                                    & Transgender Male      & 1  & 1  & 1  & 0 & 1 & 1 & 1 & 1 \\
                                    & Transgender Female    & 1  & 1  & 1  & 1 & 1 & 1 & 1 & 1 \\
                                    & Non-binary            & 1  & 1  & 1  & 0 & 1 & 1 & 1 & 1 \\
                                    & Genderqueer           & 1  & 1  & 1  & 0 & 1 & 1 & 1 & 1 \\
                                    & Intersex              & 1  & 1  & 1  & 0 & 1 & 1 & 1 & 1 \\ \midrule
\multirow{8}{*}{Sexual Orientation} & Heterosexual/straight & 1  & 1  & 1  & 1 & 1 & 1 & 1 & 1 \\
                                    & Gay                   & 1  & 1  & 1  & 0 & 1 & 1 & 1 & 1 \\
                                    & Lesbian               & 1  & 1  & 1  & 1 & 1 & 1 & 1 & 1 \\
                                    & Bisexual              & 1  & 1  & 1  & 1 & 1 & 1 & 1 & 1 \\
                                    & Queer                 & 1  & 1  & 1  & 1 & 1 & 1 & 1 & 1 \\
                                    & Questioning           & 1  & 1  & 1  & 1 & 1 & 1 & 1 & 1 \\
                                    & Asexual               & 1  & 1  & 1  & 1 & 1 & 1 & 1 & 1 \\
                                    & Pansexual             & 1  & 1  & 1  & 1 & 1 & 1 & 1 & 1 \\ \midrule
\multirow{4}{*}{Other}              & Children              & 0  & 1  & 1  & 1 & 1 & 1 & 0 & 0 \\
                                    & Citizen               & 1  & 1  & 1  & 1 & 1 & 1 & 1 & 1 \\
                                    & Immigrant             & 1  & 1  & 0  & 1 & 1 & 1 & 1 & 1 \\
                                    & HT Victim             & 1  & 1  & 1  & 1 & 1 & 1 & 1 & 1 \\ \bottomrule
\end{tabular}
\label{t:Demographics_Organization}
\end{table}

\begin{table}[]
\centering
\caption{List of service-intensity pairs used in our model and the short description of each service-intensity pair}
\begin{tabular}{c p{3cm} l p{9cm}}
\hline
$i \in I$ & Service Category               & Intensity Level & Service Description \\ 
\hline
1         & Bed                            & Single          & Crisis emergency or transitional living bed                     \\ 
\hline 
2         & \multirow{3}{3cm}{Mental Health} & Low             & Information about MH and/or stress management                    \\
3         &                                & Medium          & Talk to a counselor once a week                    \\
4         &                                & High            & Receive medication to help manage feelings                    \\ 
\hline
5         & \multirow{3}{3cm}{Medical and Dental Care} & Low   &  Medical and/or a dental check up                    \\
6         &                                & Medium          &  HIV, STI or pregnancy testing                     \\
7         &                                & High            &  Sick/injured/needs surgery                   \\ 
\hline
8         & \multirow{3}{3cm}{Substance Abuse and Alcohol Treatment} & Low             &   Information about drugs and alcohol abuse                  \\
9         &                                & Medium         & Speak to a drug/alcohol abuse counselor                     \\
10        &                                & High           & Admitted to a treatment program                    \\ 
\hline
11        & \multirow{3}{3cm}{Crisis and 24- Hour Response Services } & Low                 &  Talk to someone about stress                   \\
12        &                                & Medium          & Find somewhere to stay for this week                     \\
13        &                                & High            & Safety planning and immediate help avoid self harm                      \\ 
\hline
14        & \multirow{3}{3cm}{Long Term Support Housing }     & Low                 &  Looking for an apartment or applying for housing   \\
15        &                                & Medium             &   Housing for at least a year                \\
16        &                                & High               & Long-term place to stay in the next couple of months                       \\ 
\hline
17        & \multirow{3}{3cm}{Legal Assistance}  &Low        &Counseling to discuss rights, changing name               \\
18        &                            & Medium         & Help getting back public benefits                     \\
19        &                                & High          & Be defended in court, seek legal immigration status                    \\ 
\hline
20        & \multirow{3}{3cm}{Service Coordination}  & Low        & Advocacy for training and educational programs                     \\
21        &                                & Medium         & Advocacy for public assistance                      \\
22        &                                & High           & Advocacy for shelter and housing                      \\ 
\hline
23        & \multirow{3}{3cm}{Practical Assistance }  & Low        & Help with food, clothing or personal items 1-2/month                      \\
24        &                                & Medium         & Help with food, clothing or personal items 1-2/week                       \\
25        &                                & High           & Help with food, clothing and personal items 2+/week                      \\ 
\hline
26        & \multirow{3}{3cm}{Financial assistance }  & Low        & Financial assistance 0 -2 times per month                      \\
27        &                                & Medium         & Financial assistance 1-2 times per week                      \\
28        &                                & High           & Financial assistance 3-5 times per week                       \\ 
\hline
29         & \multirow{3}{3cm}{Life Skills} & Low             & Help learn manage responsibilities                       \\
30         &                                & Medium          & Help learn cooking, technology, staying safe                    \\
31         &                                & High            & Help learn budgeting, opening a bank account                    \\ 
\hline
32         & \multirow{3}{3cm}{Employment Assistance } & Low   &  Help learn communicate professionally   \\
33         &                                & Medium          &  Help with job training and placement  \\
34         &                                & High            & Help with resume and job searching strategies                    \\ 
\hline
35         & \multirow{3}{3cm}{Educational Assistance } & Low             &  Help signing up for GED/vocational training                  \\
36        &                                & Medium         & Help filling out financial assistance forms for college                     \\
37        &                                & High           & Help paying for GED/vocational training/college \\ 
\hline
38        & \multirow{3}{3cm}{Childcare or Parenting Help} & Low                 &  Parenting classes or coaching                     \\
39        &                                & Medium          & Occasional childcare                       \\
40        &                                & High            & Full time childcare                        \\ 
\hline
\end{tabular}
\label{t:serviceintesity}
\end{table}

\begin{table}[]
\centering
\caption{RHY organization service profiles used in modeling ($\sigma_s$), where 1 indicates the RHY organization provides a particular service-intensity}
\begin{tabular}{ccllllllll}
\hline
\multicolumn{2}{c}{\multirow{2}{*}{Service Intensity Pairs}} & \multicolumn{8}{c}{Organization} \\ \cline{3-10} 
\multicolumn{2}{c}{}                                         & 1  & 2  & 3  & 4 & 5 & 6 & 7 & 8 \\ \hline
Bed                                            & Single      & 1  & 1  & 1  & 1 & 1 & 1 & 1 &1\\
Number of beds & Single & 36 & 80 & 24 & 18 & 12 & 24 &16 & 59 \\ \hline
\multirow{3}{*}{Mental Health}                 & L           & 0  & 1  & 1  & 1 & 1 & 1 & 1 & 1 \\
                                               & M           & 0  & 1  & 1  & 1 & 1 & 1 & 1 & 1 \\
                                               & H           & 0  & 1  & 0  & 1 & 1 & 1 & 1 & 1 \\
\hline
\multirow{3}{*}{Medical}                       & L           & 0  & 0  & 0  & 0 & 0 & 1 & 1 & 1 \\
                                               & M           & 0  & 0  & 0  & 0 & 0 & 0 & 0 & 0 \\
                                               & H           & 0  & 0  & 0  & 0 & 0 & 0 & 0 & 0 \\
\hline                      \multirow{3}{*}{Substance abuse}               & L           & 0  & 1  & 1  & 1 & 1 & 1 & 1 & 1 \\
                                               & M           & 0  & 1  & 1  & 1 & 1 & 0 & 0 & 1 \\
                                               & H           & 0  & 1  & 0  & 0 & 0 & 0 & 0 & 1 \\
\hline
\multirow{3}{*}{Crisis 24 hour services}       & L           & 1  & 1  & 1  & 1 & 1 & 1 & 1 & 1 \\
                                               & M           & 1  & 1  & 1  & 1 & 1 & 1 & 1 & 1 \\
                                               & H           & 1  & 1  & 0  & 1 & 0 & 0 & 0 & 1 \\
\hline
\multirow{3}{*}{Long term housing}             & L           & 1  & 1  & 1  & 1 & 1 & 1 & 1 & 1 \\
                                               & M           & 1  & 1  & 1  & 1 & 1 & 1 & 1 & 1 \\
                                               & H           & 1  & 1  & 1  & 1 & 1 & 1 & 1 & 1 \\
\hline
\multirow{3}{*}{Legal}                         & L           & 0  & 1  & 0  & 1 & 1 & 1 & 1 & 1 \\
                                               & M           & 0  & 1  & 0  & 1 & 0 & 0 & 0 & 1 \\
                                               & H           & 0  & 1  & 0  & 1 & 0 & 0 & 0 & 0 \\
\hline
\multirow{3}{*}{Service Coordination}          & L           & 1  & 1  & 1  & 1 & 1 & 1 & 1 & 1 \\
                                               & M           & 1  & 1  & 1  & 1 & 1 & 1 & 1 & 1 \\
                                               & H           & 1  & 1  & 1  & 1 & 1 & 1 & 1 & 1 \\
\hline
\multirow{3}{*}{Practical}                     & L           & 1  & 1  & 1  & 1 & 1 & 1 & 1 & 1 \\
                                               & M           & 1  & 1  & 1  & 1 & 1 & 1 & 1 & 1 \\
                                               & H           & 1  & 1  & 1  & 1 & 1 & 1 & 1 & 1 \\
\hline
\multirow{3}{*}{Financial}                     & L           & 1  & 1  & 1  & 0 & 0 & 0 & 0 & 0 \\
                                               & M           & 1  & 1  & 1  & 0 & 0 & 0 & 0 & 0 \\
                                               & H           & 1  & 1  & 1  & 0 & 0 & 0 & 0 & 0 \\
\hline
\multirow{3}{*}{Life Skills}                   & L           & 1  & 1  & 1  & 1 & 1 & 1 & 1 & 1 \\
                                               & M           & 1  & 1  & 1  & 1 & 1 & 1 & 1 & 1 \\
                                               & H           & 1  & 1  & 1  & 1 & 1 & 1 & 1 & 1 \\
\hline
\multirow{3}{*}{Employment}                    & L           & 0  & 1  & 1  & 1 & 1 & 1 & 1 & 1 \\
                                               & M           & 0  & 1  & 0  & 1 & 1 & 1 & 1 & 1 \\
                                               & H           & 0  & 1  & 1  & 1 & 1 & 1 & 1 & 1 \\
\hline
\multirow{3}{*}{Education}                     & L           & 0  & 1  & 1  & 1 & 1 & 1 & 1 & 1 \\
                                               & M           & 0  & 1  & 0  & 1 & 1 & 1 & 1 & 1 \\
                                               & H           & 0  & 1  & 0  & 1 & 1 & 1 & 1 & 1 \\
\hline
\multirow{3}{*}{Childcare}                     & L           & 0  & 1  & 0  & 1 & 1 & 1 & 1 & 1 \\
                                               & M           & 0  & 1  & 0  & 1 & 0 & 0 & 0 & 0 \\
                                               & H           & 0  & 1  & 0  & 1 & 0 & 0 & 0 & 0
                                               \\
\hline
\end{tabular}
\label{t:Services_Organization}
\end{table}

\begin{table}[]
\centering
\caption{Percentages from the NYC Department of Youth and Community Development RHY Service 2021 report \citep{NYC-gov} and information from the National Network for Youth \citep{Monahan-2021} are used to create youth demographics profiles.}
\begin{tabular}{lcc}
\toprule
\begin{tabular}[l]{@{}c@{}}Demographic \\ Attributes\end{tabular}             & \multicolumn{1}{c}{Categories} &  \begin{tabular}[c]{@{}c@{}}Transitional \\ Living \%\end{tabular} \\ \midrule
\multirow{3}{*}{Age} & 16-17  & 9.8 \\
                     & 18-20 & 81.7 \\
                     & 21+   & 4.7 \\ \midrule
\multirow{5}{*}{Gender} & Male & 46.3 \\
                        & Female & 47.7\\
                        & Non-binary & 0.9\\
                        & Gender non-conforming& 1 \\
                      & Not sure & 0\\ \midrule
\multicolumn{2}{l}{Transgender} & 5 \\ \midrule
\multirow{8}{*}{\begin{tabular}[c]{@{}c@{}}Sexual \\ Orientation\end{tabular}} 
& Heterosexual   & 71 \\
& Gay            & 5  \\
& Lesbian        & 3  \\
& Queer          & 1 \\
& Bisexual       & 11 \\
& Asexual        & 0  \\
& Pansexual      & 2  \\
& Questioning, Not sure     & 0 \\ \midrule
Parenting Status                               & Have children     & 4  \\ \midrule
HT Victim*      & Yes   & 19-40 \\ \bottomrule
\end{tabular}
\label{t:demographics}
\end{table}
	
\begin{table}[]
\centering
\caption{The percentage of youth requesting each service, informed by interviews with stakeholders and RHY organizations.}
\begin{tabular}{lc}
\toprule
Needs of youth          & \% of youth requesting the service \\ \midrule
Bed                     & 100.0\%                            \\
Mental Health           & 16.6\%                             \\
Medical                 & 37.6\%                             \\
Substance abuse         & 6.5\%                              \\
Crisis 24 hour services & 62.7\%                             \\
Long term housing       & 56.0\%                              \\
Legal                   & 34.1\%                             \\
Service coordination    & 58.4\%                             \\
Practical               & 29.0\%                             \\
Financial               & 30.9\%                             \\
Life Skills             & 37.6\%                             \\
Employment              & 29.0\%                             \\
Education               & 62.8\%                             \\
Childcare               & 2.5\%                              \\ \bottomrule
\end{tabular}
\label{t:needs}
\end{table}

\begin{table}[]
    \centering
    \caption{The time window distributions for the services provided to youth used in our simulated data.}
    \begin{tabular}{lcc}
    \toprule
\multicolumn{3}{c}{$l_{y,i}  \thicksim Uniform(1,180)$}                \\ 
\toprule
Service Intensity Pair  & $a_{y,i}$ & $b_{y,i}$ \\
\toprule
Bed                     & $l_{y,i}$ & $a_{y,i}+ TRIA(1,2,4)$  \\
Mental Health           & $l_{y,i}$ & $a_{y,i}+ TRIA(1,2,4)$  \\
Medical                 & $l_{y,i}$ & $a_{y,i}+ TRIA(1,2,4)$  \\
Substance abuse         & $l_{y,i}$ & $a_{y,i}+ TRIA(2,3,7)$  \\
Crisis 24 hour services & $l_{y,i}$ & $a_{y,i}+ TRIA(2,3,7)$  \\
Long term housing       & $l_{y,i}$ & $a_{y,i}+ TRIA(2,3,14)$ \\
Legal                   & $l_{y,i}$ & $a_{y,i}+ TRIA(2,3,14)$ \\
Service coordination    & $l_{y,i}$ & $a_{y,i}+ TRIA(2,3,7)$  \\
Practical               & $l_{y,i}$ & $a_{y,i}+ TRIA(2,3,7)$  \\
Financial               & $l_{y,i}$ & $a_{y,i}+ TRIA(2,3,7)$  \\
Life Skills             & $l_{y,i}$ & $a_{y,i}+ TRIA(2,3,7)$  \\
Employment              & $l_{y,i}$ & $a_{y,i}+ TRIA(2,3,14)$ \\
Education               & $l_{y,i}$ & $a_{y,i}+ TRIA(2,3,14)$ \\
Childcare               & $l_{y,i}$ & $a_{y,i}+ TRIA(2,3,14)$ \\ 
    \bottomrule
    \end{tabular}
\label{t:time_windows} 
\end{table}

\end{document}